\documentclass[letterpaper, 10 pt, conference]{ieeeconf}  

\IEEEoverridecommandlockouts                              

\usepackage{comment}

\newcommand{\mydeg}[0]{$^{\circ}$}
\title{\LARGE \bf
Learning and Optimization with Bayesian Hybrid Models
}

\usepackage{color}

\usepackage{amsmath}
\usepackage{amsfonts}
\usepackage{graphicx}
\DeclareMathOperator{\E}{\mathbb{E}}
\newcommand{\D}{\mathcal{D}}
\newcommand{\vc}[1]{\mathbf{#1}}
\newcommand{\Real}{\rm I\!R}
\newcommand{\etal}[0]{et al.}

\author{Elvis A. Eugene$^{1 \dagger}$, Xian Gao$^{1 \dagger}$, and Alexander W. Dowling$^{1 *}$
\thanks{$^{1}$Department of Chemical and Biomolecular Engineering, University of Notre Dame, Notre Dame, IN 46556 USA. $^{\dagger}$ Contributed equally.
        {$^{*}$ Corresponding author: \tt\small adowling@nd.edu}}%
}

\begin{document}

\maketitle
\thispagestyle{empty}
\pagestyle{empty}

\begin{abstract}

Bayesian hybrid models fuse physics-based insights with machine learning constructs to correct for systematic bias. In this paper, we compare Bayesian hybrid models against physics-based glass-box and Gaussian process black-box surrogate models. We consider ballistic firing as an illustrative case study for a Bayesian decision-making workflow. First, Bayesian calibration is performed to estimate model parameters. We then use the posterior distribution from Bayesian analysis to compute optimal firing conditions to hit a target via a single-stage stochastic program. The case study demonstrates the ability of Bayesian hybrid models to overcome systematic bias from missing physics with less data than the pure machine learning approach. Ultimately, we argue Bayesian hybrid models are emerging paradigm for data-informed decision-making under parametric and epistemic uncertainty.

\end{abstract}

\section{INTRODUCTION}

Both \emph{model-based} and \emph{data-driven} paradigms are regularly employed for control and decision-making  \cite{Hou2013,Tidriri2016,Rudy2017,Solomatine2008}. So-called glass-box (a.k.a., white-box) models \cite{eason2016trust,eason2018advanced} are created from scientific theories, expert knowledge, and human intuition. They are regularly used to test scientific hypotheses and develop deeper fundamental understandings. When all dominant phenomena are included, glass-box models offer superior predictive power for extrapolation. Yet glass-box models are tedious to build and validate. With the availability of massive-scale data, black-box surrogate models \cite{beykal2018global,caballero2008algorithm,carpio2018kriging,fahmi2012process,jones2018superstructure,palmer2002optimization,yang2019optimization,wang2018constrained} are growing in popularity. Statistical machine learning paradigms are preferred over scientific theories for their ability to automate model (re)learning and extract undiscovered trends. At the intersection of these two paradigms are hybrid (a.k.a., grey-box) models which augments physics-based models with data-driven constructs to learn missing phenomena.

All mathematical models are burdened by uncertainty. \emph{Aleatoric} or statistical uncertainty corresponds to random phenomena, such as variability between experiments or observation noise. In contrast, \emph{epistemic} uncertainty corresponds to systematic bias or model inadequacy, such as omitting control variables in an experiment or missing physical phenomena in a model. Robust and stochastic model predictive control are well-establish paradigms to accommodate parameter uncertainty \cite{mesbah2016stochastic}. Closely related, model-based design of experiments provides a principled approach to gather information to reduce parametric uncertainty, with recent work considering approximate models \cite{franceschini2008model,quaglio2018constrained}. From a data-driven perspective, Bayesian optimization \cite{shahriari2015taking} uses probabilistic surrogate models to manage exploration (reducing model uncertainty) and exploitation (improving the objective). In this paper, we explore Kennedy-O'Hagan hybrid models as an alternative to pure glass-box and pure data-driven surrogate models.

The paper is organized as follows. In Section II, we define Bayesian hybrid models and review literature. In Section III, we demonstrate the advantages of hybrid models in illustrative case study for Bayesian calibration and optimization with epistemic uncertainty. Finally, in Section IV, we summarize our conclusions and propose future work.

\section{BAYESIAN HYBRID MODELS}

Statisticians Kennedy and O'Hagan \cite{Kennedy2001} proposed a Bayesian framework to simultaneously quantify both aleatoric and epistemic uncertainty:
\begin{equation} \label{eq:Kennedy}
	\underbrace{y_i}_{\text{observation}} = \underbrace{\eta(\mathbf{x}_i, \mathbf{\theta})}_{\text{glass-box~model}} + \underbrace{\delta(\mathbf{x}_i; \phi)}_{\text{discrepancy~model}} + \underbrace{\epsilon_i}_{\text{observation~error}}
\end{equation}

\noindent Here, observation $y_i$ is modeled with three components: a glass-box model $\eta(\cdot,\cdot)$ with physically meaningful global parameters $\mathbf{\theta}$ and inputs $\mathbf{x}_i$ for experiment $i$; a stochastic discrepancy function $\delta(\cdot)$ to counteract systematic bias in the glass-box model; and observation error $\epsilon_i$ such as uncorrelated white noise, i.e., $\epsilon_i \sim \mathcal{N}(0,\sigma^2\mathbf{I})$. Outputs of the Gaussian Process (GP) discrepancy function $\delta(\cdot) \sim \mathcal{GP}(\mathbf{\mu}(\cdot), k(\cdot,\cdot))$ follow a conditional normal distribution with the mean and covariance fully specified by $\mathbf{\mu}(\cdot)$, the mean function, and $k(\cdot,\cdot)$, the kernel function, which include hyperparameters $\phi$ \cite{Bishop2006,murphy2012machine}.
\begin{equation}
\label{eq:Bayes}
	\underbrace{p(\omega | \D)}_{\text{posterior}} = \frac{\overbrace{p(\D | \mathbf{\omega})}^{\text{likelihood}} ~ \overbrace{p(\omega)}^{\text{prior}}}{\underbrace{\int p(\D | \omega) p(\omega) d \mathbf{\omega}}_{\text{evidence}}}
\end{equation}

For Bayesian model calibration, parameters $\omega = [\theta, \phi, \sigma]$ are interpreted as random variables. One seeks to express their belief about the true value of $\omega$ using a joint probability distribution. The standard workflow is to encode one's current belief in a prior probability distribution, then observe data $\D$, and finally apply Bayes rule, Eq.~\eqref{eq:Bayes}, to obtain the posterior probability distribution $p(\omega | \D)$ \cite{Higdon2004}. For Eq.~\eqref{eq:Kennedy} and unobserved input conditions $\mathbf{x}_i$, the posterior distribution $p(\theta, \phi, \sigma | \D)$ can be propagated through $\eta(\mathbf{x}_i,\theta)$, $\delta(\mathbf{x}_i;\phi)$ and $\epsilon_i$ to obtain the posterior predictive distribution $p(y^* | \mathbf{y})$ for prediction $y_i$. Thus Kennedy and O'Hagan \cite{Kennedy2001} provide a fully Bayesian framework to calibrate global model parameters $\theta$, quantify model inadequacy (i.e., systematic bias), and make probabilistic predictions. This is an extremely powerful and flexible approach that has been used in diverse fields such as computational fluid dynamics \cite{tagade2013gaussian}, electrochemical energy storage \cite{tagade2016bayesian}, and CO$_2$ capture \cite{kalyanaraman2015bayesian,kalyanaraman2016uncertainty,Bhat2017,Li2017}.

\section{ILLUSTRATIVE COMPUTATIONAL EXAMPLE}

In this paper, we explore the performance of glass-box, data-driven, and hybrid models for decision-making. Specifically, we consider an illustrative problem of firing a ballistic to hit a desired target. We demonstrate Bayesian calibration and stochastic programming approaches for all three model classes using several training sets.

\subsection{Mathematical Models}

We start by describing the three mathematical models for the illustrative example. To avoid confusion with observations $y_i$ in Eq.~\eqref{eq:Kennedy}, we use $x$ and $z$ to represent horizontal and vertical positions, respectively. We consider a ballistic with mass $m$ (kg) fired at the origin ($x=0$, $z=0$) at angle $\psi$ (\mydeg{}) and initial velocity $v_0$ (m/s). \\

\subsubsection{Glass-box, full physics model}

We start with four coupled differential equations and initial conditions to describe the ballistic trajectory. Here $v_x$ and $v_z$ correspond with the horizontal and vertical velocities (m/s), respectively. 


\begin{equation} \label{tru1}
m \frac{dv_z}{dt} =  - mg - C_D v_z |v_z|,
\end{equation}
\begin{equation} \label{tru2}
m \frac{dv_x}{dt} =  - C_D v_x^2,
\end{equation}
\begin{equation} \label{tru7}
	\frac{d x}{dt} = v_x, \quad \frac{d z}{dt} = v_z,
\end{equation}
\begin{equation} \label{eq:initial_conditions}
\begin{split}
	v_x(0) = v_0 \cos(\psi), & \quad v_z(0) = v_0 \sin(\psi), \\ x(0) = 0, & \quad z(0) = 0 .
\end{split}
\end{equation}
The full physics glass-box model, Eqs.~\eqref{tru1} - \eqref{eq:initial_conditions}, describes projectile motion while accounting for air-resistance effects. Model parameters include $g$, acceleration due to gravity (m/s$^2$), and $C_D$, the coefficient of drag (kg/m). To account for the sign change in vertical velocity (i.e., $v_z$ decreases after being shot, $v_z$ reaches zero at the projectile's peak height, and finally $v_z$ increases during the projectile's downward motion), Eq.~\eqref{tru1} is solved in two time domains: $t \in [0, t_p)$ (upward motion) and $t \in [t_p, t_f]$ (downward motion). The absolute value in Eq.~\eqref{tru1} ensures that the drag is always in the opposite direction of the vertical velocity vector. For the first time domain (upward motion), Eq.~\eqref{tru1} becomes:
\begin{equation} \label{tru4}
m \frac{dv_z}{dt} =  - mg - C_D v_z^2,
\end{equation}
\begin{equation} \label{eq:initial_conditions2}
\begin{split}
v_z(0) = v_0 \sin(\psi), & \quad v_z(t_p) = 0, 
\end{split}
\end{equation}
The above boundary value problem is analytically solved for the time $t_p$ when the projectile reaches peak height and $v_z(t_p) = 0$. The analytic expression for $t_p$ is given in Eq.~\eqref{tru5}.
\begin{equation} \label{tru5}
t_p = \sqrt{\frac{m}{C_Dg}} \arctan \left(v_z (0) \sqrt{\frac{C_D}{mg}} \right)
\end{equation}
Next, we solve Eq.~\eqref{pkht1} to determine the peak height $z(t_p)$:
\begin{equation} \label{pkht1}
\begin{split}
\frac{d z}{dt} = v_z \\v_z(0) = v_0 \sin(\psi), & \quad v_z(t_p) = 0
\end{split}
\end{equation}
\begin{multline} \label{pkht2}
z(t_p) = \\ \frac{m}{C_D} \left( \vphantom{\left|\sqrt{\frac{}{}}\right|} \right.
\left. \ln \left( \left| \cos \left( \sqrt{\frac{C_Dg}{m}} t 
- \arctan \left( \sqrt{\frac{C_D}{mg}} v_z (0) \right) \right) \right| \right) 
\vphantom{\left|\sqrt{\frac{}{}}\right|} \right. \\
\left. - \ln \left( \left| \cos \left( \arctan \left( \sqrt{ \frac{C_D}{mg}} v_z(0) \right)\right) \right| \right) \right)
\end{multline}
The second time domain (downward motion) is described by Eqs.~\eqref{tru6} - \eqref{tru6a} and boundary conditions Eqs.~\eqref{pkht2} and \eqref{tru6ic}.
\begin{equation} \label{tru6}
m \frac{dv_z}{dt} =  - mg + C_D v_z^2
\end{equation}
\begin{equation} \label{tru6a}
\frac{d z}{dt} = v_z 
\end{equation}
\begin{equation} \label{tru6ic}
v_z(t_p) = 0,  \quad v_z(t_f) = v_f, \quad z(t_f) = 0
\end{equation}
The analytic solution for this system is given by Eq. \eqref{tru8}.
\begin{multline}\label{tru8}
z(t_f) = \sqrt{\frac{mg}{C_D}}\left(t_p - t_f + \left(\ln(2) \vphantom{\sqrt{\frac{m}{C_Dg}}} \right. \right. \\
\left. \left. - \sqrt{\frac{m}{C_Dg}} 
\ln \left(\exp \left(2 \sqrt{\frac{C_Dg}{m}} (t_p - t_f) \right) + 1 \right) \right) \right) + z(t_p)
\end{multline}

To compute $t_f$, we set $z(t_f) = 0$ and numerically solve Eq. \eqref{tru8}. Finally, we use $t_f$ in Eq.~\eqref{tru3} to compute the total horizontal distance traveled. Our analysis assumes no change in elevation.
\begin{equation} \label{tru3}
\hat{y}(t_f) = \frac{m}{C_D} \ln \left( C_D v_0 \cos{\psi}) t_f + m\right) 
\end{equation}
\begin{equation} \label{eq:train}
y (t_f) = \hat{y}(t_f) + \epsilon, \quad \epsilon \sim \mathcal{N}(0,\sigma^2 \vc I)
\end{equation}
In this paper, we use Eq.~\eqref{eq:train}, which includes observation noise, in place of physical experiments. 

We now describe three predictive surrogate models with inputs $\mathbf{x}_i = [v_{0i}, \psi_i]$ where $v_{0i}$ is the initial velocity and $\psi_i$ is the launching angle for experiment $i$. The input dimensionality is $D=2$.
\subsubsection{Glass-box, simple physics surrogate model} Often in engineering applications, the full physics model is not known, computationally intractable, or otherwise cumbersome to work with. Instead, a simplified surrogate model is often used. For the ballistics problem, we consider a parabolic trajectory that neglects air resistance:
\begin{equation} \label{simp}
\begin{aligned}
 y_i &= \eta(\mathbf{x}_i, \mathbf{\theta})  + \epsilon_i\\
 \eta(\mathbf{x}_i, \mathbf{\theta}) &= \frac{2(v_{0i})^2}{g} \sin(\psi_i) \cos(\psi_i)
\end{aligned}
\end{equation}
The impact location $y$ (horizontal distance traveled) given by Eq.~\eqref{simp} does not depend on $m$, the projectile mass. \\

\subsubsection{Black-box, Gaussian Process (GP) surrogate model}
Another approach is to construct a surrogate model completely from data, without incorporating (significant) knowledge or physical intuition. Here, we will consider $y_i \sim \mathcal{GP}(\mathbf{\mu}(\cdot), k(\cdot,\cdot))$, i.e., a GP model with observation error:
\begin{equation} \label{eq:pure_gp}
y_i= \delta(\mathbf{x}_i; \vc {\phi}) + \epsilon_i.
\end{equation}

We use a Gaussian distribution with a zero mean and a covariance ($\vc K$) constructed by the kernel function ($k$) as a prior distribution for $\delta$:
\begin{equation} \label{gaussian_prior}
p(\vc \delta) = \mathcal{N}(\vc \delta | \vc 0,\vc K)
\end{equation}
where $\vc K_{mn} = k(\vc x_m, \vc x_n)$. We use the radial-based function (RBF) with automatic relevance determination (ARD) \cite{williams2006gaussian}:
\begin{equation}
k(\vc x_m, \vc x_n) = \sigma_f^2 \exp[-\frac{1}{2} \sum_{i=1}^{D} \frac{(x_{mi} - x_{ni})^2}{l_i^2}]
\end{equation}
where $\sigma_f^2$ is the kernel variance and $l_i$ is the lengthscale. Also a Gaussian likelihood is assumed for the $N$ training data points:
\begin{equation} \label{gaussian_likelihood}
p(\vc y |\vc \delta) = \mathcal{N}(\vc y | \vc \delta, \sigma^{2} \vc I_N)
\end{equation}
where $\sigma^2$ is the variance of the observation noise. Using the properties of linear Gaussian models, one can derive the marginal likelihood function \cite{Bishop2006}:

\begin{equation} \label{marginal_likelihood}
p(\vc y) =\int p(\vc y |\delta) p(\vc \delta) d\vc \delta =  \mathcal{N}(\vc y | \vc 0, \vc C)
\end{equation}
where $\vc C = \vc K + \sigma^{2}\vc I_N$. Conditioning on the training data, we then derive the predictive distribution:

\begin{equation} \label{predictive_distribution}
p(y^{*}|\vc y) = \mathcal{N}(y^{*}| m_{GP}(\vc x^{*}), \sigma_{GP}^2(\vc x^{*}))
\end{equation}
Here, $y^*$ is the predicted observation for a new input $x^*$ that is not contained in the training data. Notice $y^{*}$ is a random scalar and its predictive distribution is Gaussian with predictive mean $m_{GP}(\vc x^{*})$ and predictive variance $\sigma_{GP}^2(\vc x^{*})$ functions: 

\begin{equation} \label{predictive_functions}
\begin{aligned}
m_{GP}(\vc x^{*}) &= \vc k^T \vc C^{-1} \vc y\\
\sigma^2_{GP}(\vc x^{*}) &= c - \vc k^T \vc C^{-1} \vc k
\end{aligned}
\end{equation}
These functions depend on $\vc x^{*}$ through $c = k(\vc x^*,\vc x^{*})$ and $\vc k = [k(\vc x_1,\vc x^{*}),...,k(\vc x_N,\vc x^{*})]^T$.\\

\subsubsection{Hybrid surrogate model} As we will see later, there are advantages and disadvantages to both physics-based and data-driven surrogate models. We now consider a Kennedy-O'Hagan hybrid model:

\begin{equation} \label{hybrid_model}
\begin{aligned} 
y_i= \eta(\mathbf{x}_i, \mathbf{\theta}) + \delta(\mathbf{x}_i; \phi)+ \epsilon_i, \\
\eta(\mathbf{x}_i, \mathbf{\theta}) = \frac{2v_0^2}{g} \sin(\psi) \cos(\psi)\\
\delta(\mathbf{x}_i; \phi)+ \epsilon_i \sim \mathcal{GP}(\mu(\cdot),k(\cdot,\cdot))
\end{aligned}
\end{equation}
Here, we combine the simple physics surrogate model Eq.~\eqref{simp} and the data-driven Gaussian process Eqs. \eqref{eq:pure_gp} - \eqref{predictive_functions} to establish the Bayesian hybrid model Eq.~\eqref{hybrid_model}. We note that the GP surrogate is now trained on the difference between the observations $\mathbf{y}$ and the simple physics model outputs $\eta(\mathbf{x}, \mathbf{\theta})$.


\subsection{Training Data}

\begin{figure}[!ht]
	\centering
	\includegraphics[width=0.5\textwidth]{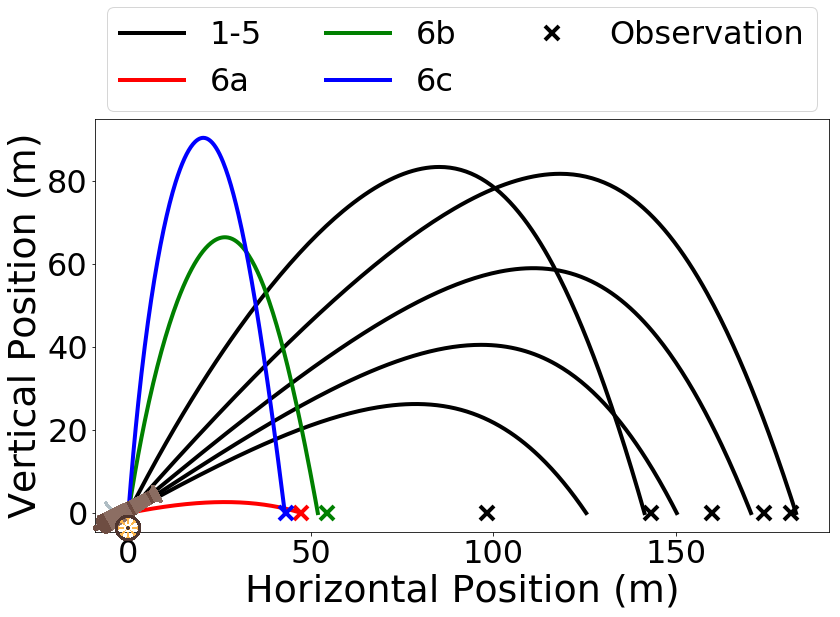}
	\caption{Trajectories for training data generated from the full physics model. Significant air resistance causes the trajectories to be not parabolic. \texttt{x} marks the observed distance traveled and includes observation error.} \label{fig:training}
\end{figure}

We assume $g = 9.8$ m/s$^2$ and $C_D/m = 0.01$ kg/m${}^2$ and simulate the full physics model, Eq. \eqref{tru1} - \eqref{tru3}, to generate training data. Each experiment consists of the initial velocity $v_0$, firing angle $\psi$, and observed impact location $y$ (horizontal distance traveled). We add Gaussian observation noise with mean zero and standard deviation 5 m to represent aleatoric uncertainty. Figure \ref{fig:training} and Table \ref{tab:train} summarizes the training data. In this paper, experiments 1-5 act as a base training dataset. We compare the impact of selecting either experiment 6a, 6b, or 6c as the sixth experiment on both Bayesian model calibration and optimization under uncertainty.

\begin{table}[!ht]
	\centering
	\caption{Training data generated from full physics model.} \label{tab:train}
	
	\begin{tabular}{c c c c}
		\textbf{Experiment} & \textbf{Angle (\mydeg)} & \textbf{Initial velocity (m/s)} & \textbf{Impact location (m)}\\
		{} & $\psi$  & $v_0$ & $y$ \\ \hline
		1 & 25 & 60 & 118.18 \\
		2 & 30 & 70 & 159.79 \\
		3 & 36 & 80 & 174.14 \\
		4 & 45 & 90 & 181.67 \\
		5 & 60 & 75 & 143.21 \\
		6a & 10 & 42 & 47.305 \\
		6b & 80 & 53 & 54.294 \\
		6c & 85 & 71 & 43.239 \\ \hline
	\end{tabular}
	
\end{table}

\subsection{Bayesian Calibration}

For each of the three postulated models, we wish to use a subset of the training data to infer the unknown model parameters. We use a Bayesian approach, i.e., Eq.~\eqref{eq:Bayes}, which takes a prior distribution and computes a posterior distribution. For all three surrogate models, Bayesian calibration is performed using Markov Chain Monte Carlo implemented in Python using the PyMC3 package \cite{salvatier2016probabilistic}.

\subsubsection{Simple physics surrogate model}

In the simple physics model, Eq.~\eqref{simp}, we seek to estimate parameter $g$. We use a uniform prior distribution:
\begin{equation}
\frac{1}{g} \sim \mathcal{U}\left(0.001,1\right)
\end{equation}

The likelihood $p(\D|\omega)$ is assumed to follow a normal distribution. This is equivalent to assuming the observation error is normally distributed.
\begin{equation}
p(\D|\omega) = \mathcal{N}\left(\eta(\mathbf{x_i},g),\tau^{-1}\right)
\end{equation}

Here, $\eta(\mathbf{x_i},g)$ is the simple model and $\tau$ is the precision, i.e., inverse of the observation error variance. We assume a gamma prior for hyperparameter $\tau$:
\begin{equation}
\tau \sim \mathcal{G}\left(0.25,2.5\right)
\end{equation}

For the simple model, only one model parameter ($g$) and one model hyperparameter ($\tau$) needs to be estimated, thus:
\begin{equation}
\omega = [g, \tau]
\end{equation}

\subsubsection{GP surrogate model} Recall for the ballistics problem, $D = 2$. The hyperparameters $\vc \phi = [\sigma_f^2, l_1, l_2, \sigma^{2}]$ includes the kernel variance $\sigma_f^{2}$, length scales $l_1$ and $l_2$ and likelihood variance $\sigma^2$. The hyperparameters are given the following prior distributions:
\begin{equation}\label{hyperparam_prior}
\begin{aligned}
\sigma_f \sim \mathcal{U} (0.1,1) \\
l_i \sim \mathcal{U} (1,50) \\
\sigma \sim \mathcal{N}(0,5)
\end{aligned}
\end{equation}

We combine these prior distributions with the marginal likelihood, Eq. \ref{marginal_likelihood}, which gives the the unnormalized posterior. We then compute a maximum a posteriori (MAP) point estimate for the four GP hyperparameters:
\begin{equation}
\phi = [\sigma_f^2, l_1, l_2, \sigma^{2}]
\end{equation}

With the estimated $\phi$ and an arbitrary input $\vc x^*$, predictive mean and variance can be calculated by Eq.~\eqref{predictive_functions}.

\subsubsection{Hybrid surrogate model} Calibrating the hybrid model requires inferring the parameters in both the glass-box and data-driven components. Thus,
\begin{equation}
\theta = g, \quad \phi = [\sigma_f^2, l_1, l_2, \sigma^{2}]
\end{equation}
To train this hybrid model, we first infer $g$ using Bayesian regression. We then compute the residual $y_i - \eta(\mathbf{x}_i, g)$. Finally, we train the GP via MAP estimation on the residual with the prior distributions in Eq.~\eqref{hyperparam_prior}. 

\subsection{Optimization under Uncertainty}

We now use the posterior distribution from Bayesian calibration for optimization under uncertainty. Specifically, we seek a velocity $v_0$ and angle $\psi$ pair to hit a target at $\bar{y}$ = 100 m, which we formulate as a single-stage stochastic program:
\begin{equation} \label{eq:cannon-OUU}
	\max_{v_0,\psi} ~~ \mathbb{E}[u(\bar{y} - y_j)], \quad u(\Delta y) = 1 - \frac{1}{100} \min{} (\Delta y, 100) 
\end{equation}

The utility function $u(\cdot)$ varies linearly from unity for a direct hit to zero for a miss of 100 meters or more. We approximate the expectation in Eq.~\eqref{eq:cannon-OUU} using 4500 samples from the posterior distribution of $\theta$. For computational simplicity, we use a MAP point estimate for hyperparameters $\phi$. For all three surrogate models, we use a Gauss-Hermite quadrature rule with 7 nodes to approximate the expected value from the GP prediction (if applicable) and observation error distributions. We perform optimization by grid search from $v_0 = $ 40 to 100 m/s and $\psi = $ 1\mydeg{} to 90\mydeg{} to facilitate visualization.


\subsection{Results}
\begin{table*}[!ht]
	\centering
	\caption{Results for calibration and optimization with nine model and training dataset combinations.} \label{res-tab}
	
	\begin{tabular}{c c c c c c c}
		\textbf{Training} & \textbf{Training} & \textbf{Calibrated} & \textbf{Optimum} & \textbf{Optimum} & \textbf{Expected} & \textbf{Observed}\\ 
		\textbf{set} & \textbf{experiments} & \textbf{model} & \textbf{angle (\mydeg)} & \textbf{velocity (m/s)} & \textbf{distance (m)} & \textbf{distance (m)}\\
		\textbf{} & \textbf{} & \textbf{} & \textbf{$\psi^{\dagger}$} & \textbf{$v_0^{\dagger}$} & \textbf{} & \textbf{$y$}\\\hline
		{} & {} & {Simple} & {40} & {60.0} & {97.0} & {141} \\ 
		{A} & {1,2,3,4,5,6a} &  {GP} & {29} & {57.5} & {137} & {120} \\ 
		{} & {} &  {Hybrid} & {22} & {57.5} & {77.7} & {113} \\ \hline 
		{} & {} &  {Simple} & {18} & {77.5} & {97.0} & {139} \\ 
		{B} & {1,2,3,4,5,6b} &  {GP} & {31} & {57.5} & {138} & {132} \\ 
		{} & {} &  {Hybrid} & {22} & {57.5} & {77.7} & {118} \\ \hline  
		{} & {} &  {Simple} & {13} & {90.0} & {96.0} & {139} \\ 
		{C} & {1,2,3,4,5,6c} &  {GP} & {26} & {60.0} & {137} & {118} \\ 
		{} & {} & {Hybrid} & {72} & {72.5} & {97.7} & {116} \\ \hline

	\end{tabular}
	
\end{table*}
We now compare the Bayesian calibration and optimization under uncertainty workflow for combinations of three training datasets (A, B, C) and three models (simple physics, GP, hybrid). For each training dataset, we perform Bayesian calibration and then compute $v_0^{\dagger}$ and $\psi^{\dagger}$ that maximize Eq.~\eqref{eq:cannon-OUU}. Each training dataset has six experiments. The first five experiments are consistent across all training sets. Table \ref{res-tab} shows how changing the 6th experiment impacts the optimization results. From Table \ref{res-tab}, we observe several trends:
\begin{enumerate}
	\item The \textbf{glass-box simple physics surrogate model} behaves poorly and consistently overshoots the target by at least 39 m each time. This is due to the fact that the simple model fails to account for all physical phenomena governing projectile motion, specifically, air-resistance effects. The maximum objective value calculated using the simple physics model to find the optimum firing angle and velocity is 0.698. Recall the objective value 1 indicates a direct hit of the target with exact certainty.
	\item The \textbf{black-box GP surrogate model} also consistently misses the target, but the maximum distance overshot is 32 m which is less than the simple model. The superior performance of the GP model can be attributed to its purely data-driven nature, with no reliance on the physics governing the process. The maximum objective value calculated when using the pure GP model is 0.818.
	\item The \textbf{hybrid surrogate model} outperforms the simple and GP models. Out of the nine model and data set combinations, the maximum overshooting is restricted to 18 m. We attribute the success of the hybrid model to its ability to learn missing physics with less data from the pure data-driven approach. The maximum objective value obtained when using the hybrid model is 0.987.
\end{enumerate}

\begin{figure}[!ht]
	\centering
	\includegraphics[width=0.45\textwidth]{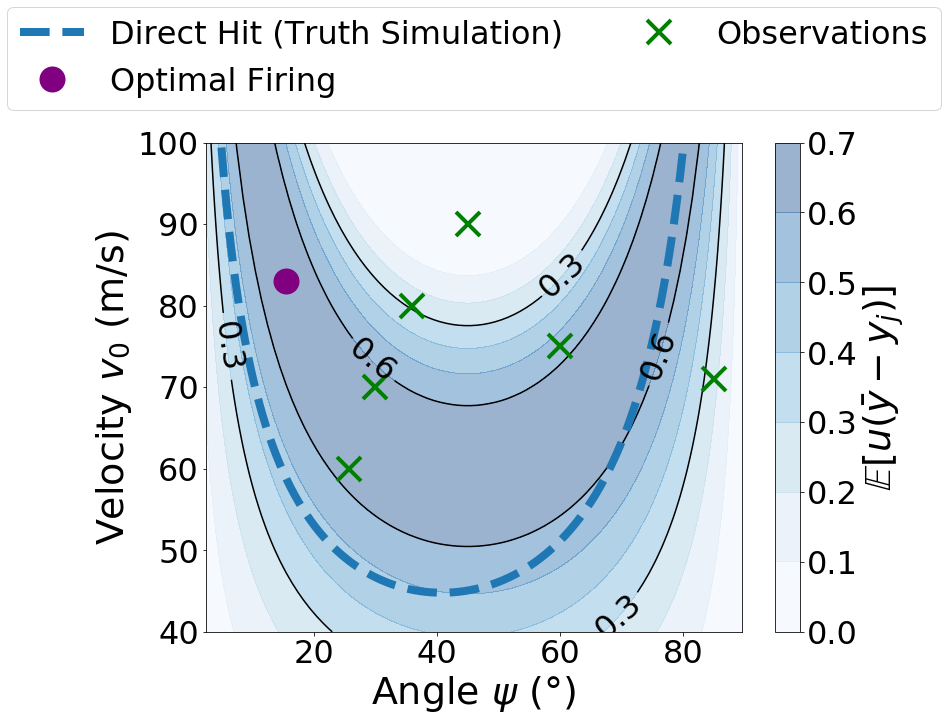}
	\caption{The wide contours of the objective function $\E[u]$ indicate high uncertainty using the \textbf{simple model}. The highest objective value achieved is 0.698 ($\E[u]=1$ indicates a direct hit with probability 1).} \label{res-simp}
\end{figure}

\begin{figure}[!ht]
	\centering
	\includegraphics[width=0.45\textwidth]{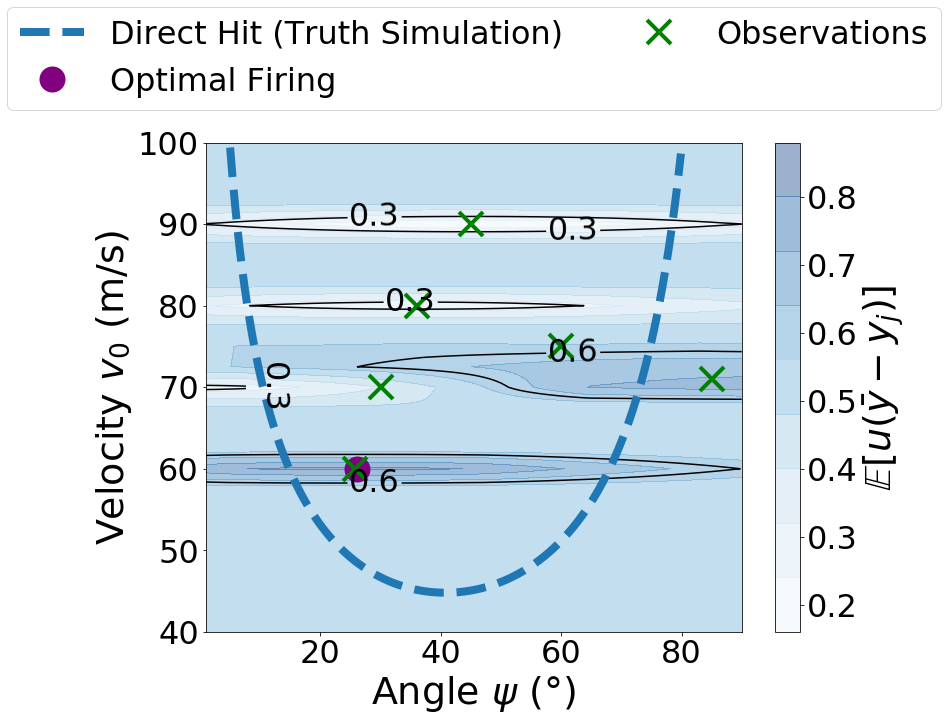}
	\caption{The narrower contours of the objective function with the \textbf{GP model} indicate lower uncertainty than the simple model. The highest objective value calculated is 0.818 ($\E[u]=1$ indicates a direct hit with probability 1). The model is insensitive to variations in firing angle indicating the need to refine the training procedure.} \label{res-gp}
\end{figure}

\begin{figure}[!ht]
	\centering
	\includegraphics[width=0.45\textwidth]{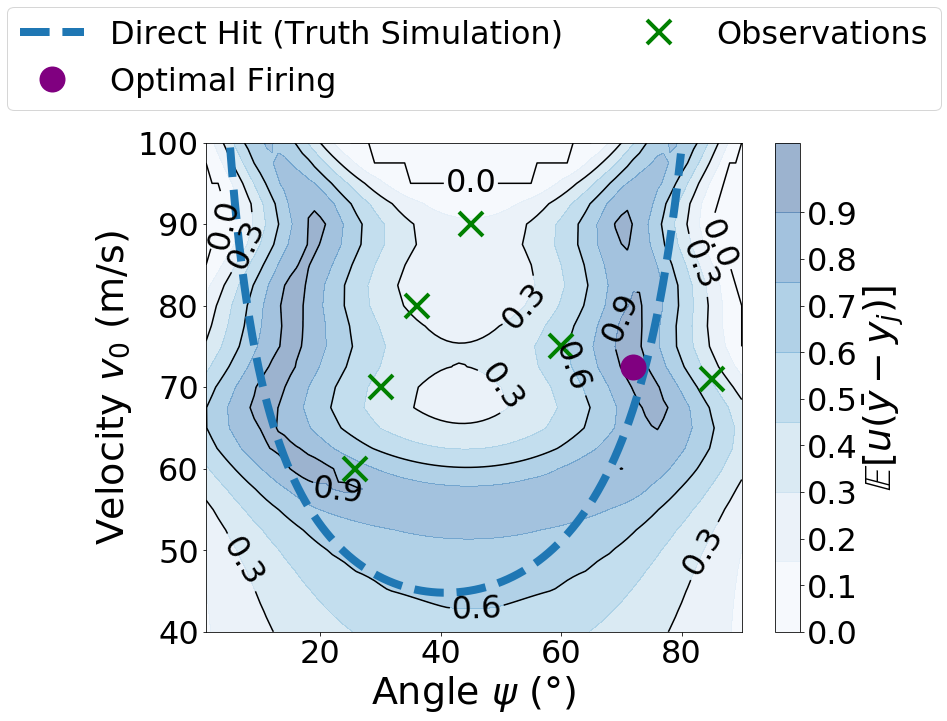}
	\caption{The \textbf{hybrid model} performs the best out of the three models with the highest objective value of 0.987 ($\E[u]=1$ indicates a direct hit with probability 1). The physics from the simple model makes the model sensitive to the firing angle and the GP accounts for the neglected air resistance in the simple model to provide the optimal firing conditions.} \label{res-hyb}
\end{figure}


Next, we examine training dataset C to further compare the models. Figures \ref{res-simp} - \ref{res-hyb} show the contours of the objective function in Eq.~\eqref{eq:cannon-OUU} for all three models. Recall, the utility function $u(\cdot)$ varies linearly between 1.0 for a direct hit and 0.0 for a miss of 100 m or more. The dashed blue line shows $\psi$, $v_0$ combinations for a direct hit with the full physics model, i.e., the truth simulation. We make several observations from these plots:
\begin{enumerate}
	\item In Figure \ref{res-simp} we see wide contours for the objective that peak below 0.7. This indicates that there is high uncertainty with the \textbf{simple model}. Moreover, the average value for $g$ in the posterior is 35.933 m/s$^2$. The estimate of $g$ is biased to offset the missing air resistance in the model.
	\item While Figure \ref{res-gp} shows higher objective values (above 0.8), it also reveals oblong contours. The \textbf{GP model} is fairly insensitive to angle $\psi$. This suggests further refinement of the kernel and/or input transformations should be explored.
	\item Figure \ref{res-hyb} illustrates the superior performance of the \textbf{hybrid model}. The glass-box component provides the basic shape and the data-driven component corrects for missing physics. As expected, the contours of the hybrid model best align with the truth simulation in regions near training data.
\end{enumerate}

\section{CONCLUSIONS}

Through the ballistics illustrative example, we compare three surrogate modeling paradigms: glass-box, data-driven, and hybrid. We establish Bayesian workflow for data-informed decision-making under uncertainty, where the posterior distribution from Bayesian inference is directly used for stochastic programming. We see how model-form uncertainty, specifically neglecting air resistance, biases both parameter estimates ($g$ is 4 times larger than reality) and predictions (always overshoots). Moreover, a pure data-driven approach using GP regression requires careful tuning and struggles with limited data. We found the Kennedy-O'Hagan hybrid models overcome limitations of both of these approaches.

More broadly, Kennedy-O'Hagan hybrid models provide a new direction to account for epistemic (i.e., model-form) uncertainty in stochastic control and decision-making. We specifically advocate for a Bayesian approach, which provides a rich statistical framework to obtain probability distributions from data that are needed for stochastic programming. As future work, we plan to expand the ballistics case study to understand how each modeling approach performs with additional data. We seek to compare sequential and simultaneous inference of parameters $\theta$ and hyperparameters $\phi$. We hypothesize the latter can reduce bias. We are also interested in other acquisition functions to manage the trade-off between exploration (reducing model uncertainty) and exploitation (direct hit) \cite{shahriari2015taking}. Finally, we plan to consider large-scale applications such as chemical reactor kinetic modeling, design, and control.

\addtolength{\textheight}{-12cm}   






\bibliographystyle{ieeetr}
\bibliography{crii,MyPublications,Diafiltration,career}

\begin{thebibliography}{10}

\bibitem{Hou2013}
Z.-S. Hou and Z.~Wang, ``From model-based control to data-driven control:
  Survey, classification and perspective,'' {\em Information Sciences},
  vol.~235, pp.~3 -- 35, 2013.
\newblock Data-based Control, Decision, Scheduling and Fault Diagnostics.

\bibitem{Tidriri2016}
K.~Tidriri, N.~Chatti, S.~Verron, and T.~Tiplica, ``Bridging data-driven and
  model-based approaches for process fault diagnosis and health monitoring: A
  review of researches and future challenges,'' {\em Annual Reviews in
  Control}, vol.~42, pp.~63 -- 81, 2016.

\bibitem{Rudy2017}
S.~H. Rudy, S.~L. Brunton, J.~L. Proctor, and J.~N. Kutz, ``Data-driven
  discovery of partial differential equations,'' {\em Science Advances},
  vol.~3, no.~4, 2017.

\bibitem{Solomatine2008}
D.~P. Solomatine and A.~Ostfeld, ``Data-driven modelling: some past experiences
  and new approaches,'' {\em Journal of Hydroinformatics}, vol.~10, no.~1,
  p.~3, 2008.

\bibitem{eason2016trust}
J.~P. Eason and L.~T. Biegler, ``A trust region filter method for glass
  box/black box optimization,'' {\em AIChE Journal}, vol.~62, no.~9,
  pp.~3124--3136, 2016.

\bibitem{eason2018advanced}
J.~P. Eason and L.~T. Biegler, ``Advanced trust region optimization strategies
  for glass box/black box models,'' {\em AIChE Journal}, vol.~64, no.~11,
  pp.~3934--3943, 2018.

\bibitem{beykal2018global}
B.~Beykal, F.~Boukouvala, C.~A. Floudas, N.~Sorek, H.~Zalavadia, and E.~Gildin,
  ``Global optimization of grey-box computational systems using surrogate
  functions and application to highly constrained oil-field operations,'' {\em
  Computers \& Chemical Engineering}, vol.~114, pp.~99--110, 2018.

\bibitem{caballero2008algorithm}
J.~A. Caballero and I.~E. Grossmann, ``An algorithm for the use of surrogate
  models in modular flowsheet optimization,'' {\em AIChE journal}, vol.~54,
  no.~10, pp.~2633--2650, 2008.

\bibitem{carpio2018kriging}
R.~R. Carpio, F.~F. Furlan, R.~C. Giordano, and A.~R. Secchi, ``A kriging-based
  approach for conjugating specific dynamic models into whole plant stationary
  simulations,'' {\em Computers \& Chemical Engineering}, vol.~119,
  pp.~190--194, 2018.

\bibitem{fahmi2012process}
I.~Fahmi and S.~Cremaschi, ``Process synthesis of biodiesel production plant
  using artificial neural networks as the surrogate models,'' {\em Computers \&
  Chemical Engineering}, vol.~46, pp.~105--123, 2012.

\bibitem{jones2018superstructure}
M.~Jones, H.~Forero-Hernandez, A.~Zubov, B.~Sarup, and G.~Sin, ``Superstructure
  optimization of oleochemical processes with surrogate models,'' in {\em
  Computer Aided Chemical Engineering}, vol.~44, pp.~277--282, Elsevier, 2018.

\bibitem{palmer2002optimization}
K.~Palmer and M.~Realff, ``Optimization and validation of steady-state
  flowsheet simulation metamodels,'' {\em Chemical Engineering Research and
  Design}, vol.~80, no.~7, pp.~773--782, 2002.

\bibitem{yang2019optimization}
S.~Yang, S.~Kiang, P.~Farzan, and M.~Ierapetritou, ``Optimization of reaction
  selectivity using {CFD}-based compartmental modeling and surrogate-based
  optimization,'' {\em Processes}, vol.~7, no.~1, p.~9, 2019.

\bibitem{wang2018constrained}
Z.~Wang and M.~Ierapetritou, ``Constrained optimization of black-box stochastic
  systems using a novel feasibility enhanced kriging-based method,'' {\em
  Computers \& Chemical Engineering}, vol.~118, pp.~210--223, 2018.

\bibitem{mesbah2016stochastic}
A.~Mesbah, ``Stochastic model predictive control: An overview and perspectives
  for future research,'' {\em IEEE Control Systems}, vol.~36, no.~6,
  pp.~30--44, 2016.

\bibitem{franceschini2008model}
G.~Franceschini and S.~Macchietto, ``Model-based design of experiments for
  parameter precision: State of the art,'' {\em Chemical Engineering Science},
  vol.~63, no.~19, pp.~4846--4872, 2008.

\bibitem{quaglio2018constrained}
M.~Quaglio, E.~S. Fraga, and F.~Galvanin, ``Constrained model-based design of
  experiments for the identification of approximated models,'' {\em
  IFAC-PapersOnLine}, vol.~51, no.~15, pp.~515--520, 2018.

\bibitem{shahriari2015taking}
B.~Shahriari, K.~Swersky, Z.~Wang, R.~P. Adams, and N.~De~Freitas, ``Taking the
  human out of the loop: A review of bayesian optimization,'' {\em Proceedings
  of the IEEE}, vol.~104, no.~1, pp.~148--175, 2015.

\bibitem{Kennedy2001}
M.~C. Kennedy and A.~O'Hagan, ``{Bayesian calibration of computer models},''
  {\em Journal of the Royal Statistical Society: Series B (Statistical
  Methodology)}, vol.~63, no.~3, pp.~425--464, 2001.

\bibitem{Bishop2006}
C.~M. Bishop, {\em Pattern Recognition and Machine Learning}.
\newblock Springer, 2006.

\bibitem{murphy2012machine}
K.~P. Murphy, {\em Machine Learning: A Probabilistic Perspective}.
\newblock MIT press, 2012.

\bibitem{Higdon2004}
D.~Higdon, M.~Kennedy, J.~C. Cavendish, J.~A. Cafeo, and R.~D. Ryne,
  ``Combining field data and computer simulations for calibration and
  prediction,'' {\em SIAM Journal on Scientific Computing}, vol.~26, no.~2,
  pp.~448--466, 2004.

\bibitem{tagade2013gaussian}
P.~M. Tagade, B.-M. Jeong, and H.-L. Choi, ``A {Gaussian} process emulator
  approach for rapid contaminant characterization with an integrated
  {multizone-CFD} model,'' {\em Building and Environment}, vol.~70,
  pp.~232--244, 2013.

\bibitem{tagade2016bayesian}
P.~Tagade, K.~S. Hariharan, S.~Basu, M.~K.~S. Verma, S.~M. Kolake, T.~Song,
  D.~Oh, T.~Yeo, and S.~Doo, ``Bayesian calibration for electrochemical thermal
  model of lithium-ion cells,'' {\em Journal of Power Sources}, vol.~320,
  pp.~296--309, 2016.

\bibitem{kalyanaraman2015bayesian}
J.~Kalyanaraman, Y.~Fan, Y.~Labreche, R.~P. Lively, Y.~Kawajiri, and M.~J.
  Realff, ``Bayesian estimation of parametric uncertainties, quantification and
  reduction using optimal design of experiments for {CO$_2$} adsorption on
  amine sorbents,'' {\em Computers \& Chemical Engineering}, vol.~81,
  pp.~376--388, 2015.

\bibitem{kalyanaraman2016uncertainty}
J.~Kalyanaraman, Y.~Kawajiri, R.~P. Lively, and M.~J. Realff, ``Uncertainty
  quantification via {Bayesian} inference using sequential {Monte Carlo}
  methods for {CO$_2$} adsorption process,'' {\em AIChE Journal}, vol.~62,
  no.~9, pp.~3352--3368, 2016.

\bibitem{Bhat2017}
K.~S. Bhat, D.~S. Mebane, P.~Mahapatra, and C.~B. Storlie, ``Upscaling
  uncertainty with dynamic discrepancy for a multi-scale carbon capture
  system,'' {\em Journal of the American Statistical Association}, vol.~112,
  no.~520, pp.~1453--1467, 2017.

\bibitem{Li2017}
K.~Li, P.~Mahapatra, K.~S. Bhat, D.~C. Miller, and D.~S. Mebane, ``Multi-scale
  modeling of an amine sorbent fluidized bed adsorber with dynamic discrepancy
  reduced modeling,'' {\em Reaction Chemistry \& Engineering}, vol.~2, no.~4,
  pp.~550--560, 2017.

\bibitem{williams2006gaussian}
C.~K. Williams and C.~E. Rasmussen, {\em Gaussian Processes for Machine
  Learning}, vol.~2.
\newblock MIT press Cambridge, MA, 2006.

\bibitem{salvatier2016probabilistic}
J.~Salvatier, T.~V. Wiecki, and C.~Fonnesbeck, ``Probabilistic programming in
  {Python} using {PyMC3},'' {\em PeerJ Computer Science}, vol.~2, p.~e55, 2016.

\end{thebibliography}

\end{document}